\documentclass[11pt]{amsart}
\usepackage{amssymb, amsmath,amsthm}

\renewcommand{\epsilon}{\varepsilon}

\newtheorem{theorem}{Theorem}

\begin{document}
\parindent0 in
\parskip 1 em
\title{Quantum Layers over Surfaces Ruled Outside a Compact Set}
\author{Christopher Lin and Zhiqin Lu}

\address{Department of
Mathematics, University of California,
Irvine, Irvine, CA 92697}

\email[Christopher Lin]{ccplin@msri.org}
\email[Zhiqin Lu]{zlu@math.uci.edu}
\begin{abstract}
In this paper, we proved the quantum layer over a surface which is ruled outside a compact set, asymptotically flat but not totally geodesic admits ground states.
\end{abstract}

\thanks{The
second author is partially supported by  NSF Career award DMS-0347033 and the
Alfred P. Sloan Research Fellowship.}

\maketitle

\section{Introduction}\label{sec1}
\quad The spectrum of the Laplacian on manifolds is a classic domain of research within geometric analysis.  One of its least developed area is spectral analysis on noncompact, non-complete manifolds.  An insteresting paper by Duclos, Exner, and Krej\v{c}i\v{r}\'{\i}k 
(\cite{dek}, 2001) demonstrated the existence of discrete spectrum, or bound states, of the Dirichlet Laplacian acting on $L^2$ functions on a particular type of noncompact, non-complete manifold which they termed the ``quantum layer''.  In particular, the existence of bound states means that at least there exists the lowest bound state - the ground state.  The main results in their paper were generalized in \cite{LL-1} and 
\cite{LL-2} by Lin and Lu to higher dimensions under more general geometric settings.  The existence of discrete spectrum for the Laplacian is a non-trivial phenomenon on noncompact manifolds, even when the manifold is complete. The results in \cite{dek}, \cite{LL-1}, and 
\cite{LL-2} are rare instances in which the discrete spectrum is clearly known to exist on noncompact, non-complete manifolds.

\quad We recall the definition of a quantum layer below, following \cite{LL-1}. 
\newtheorem{definition}{Definition} 
\begin{definition}
Let $\varSigma \hookrightarrow \mathbb{R}^{3}$ be an isometrically immersed, oriented  hypersurface.  Let $N$ and $\vec{A}$ respectively be the unit normal vector field and the second fundamental form on $\varSigma$, and define the map 
\begin{equation}\label{layerimm}
p: \varSigma \times (-a,a) \longrightarrow \mathbb{R}^{3}
\end{equation}
by $(x,u) \mapsto x + uN$, where the number $a>0$, which is called the thickness, is such that there is a constant 
$C_o$ such that $a\|\vec{A}\| < C_o <1$ on $\varSigma$.  We define a quantum layer to be the smooth manifold $\Omega = \varSigma \times (-a,a)$ with the pull-back metric $p*(ds^2_E)$.   
\end{definition} 
As one can see from the definition, a quantum layer is just a tubular neighborhood of a surface in $\mathbb{R}^3$.  The relation between $\vec{A}$ and the thickness $a$ above is simply to ensure that the map (\ref{layerimm}) is also an immersion.  
The terminology of ``quantum'' simply eludes to the fact that when suitable geometric conditions are imposed on the tubular neighborhood, ground state exists.   

\quad Before we go any further, let us establish a consensus on the notions and terminologies from functional analysis used in this paper.  A good source of this can, for example, be found in \cite{reed-simon-1}.  We will use $\Delta$ to denote the Dirichlet Laplacian throughout the paper, which is a self-adjoint exension in $L^2(\Omega)$ of the positive-definite Laplacian $\Delta = -\text{div}\circ\text{grad}$ on smooth, compactly supported functions on $\Omega$.  The Dirichlet Laplacian $\Delta$ is thus an unbounded operator with a dense domain 
$\text{Dom}(\Delta)\subset W_0^{1,2}(\Omega) \subset L^2(\Omega)$.  
The resolvent set $\rho(\Delta)$ of the operator is the set of all complex numbers such that 
$(\Delta -\lambda I)^{-1}$ exists as a bounded operator on $L^2(\Omega)$, and the spectrum $\sigma(\Delta)$ is defined as the set $\mathbb{C}\setminus\rho(\Delta)$.  The self-adjointness and positive-definiteness of $\Delta$ imply that $\sigma(\Delta)\subset [0,\infty)$.  The spectrum can be decomposed as 
$\sigma(\Delta) = \sigma_{disc}(\Delta)\sqcup\sigma_{ess}(\Delta)$, where the discrete spectrum $\sigma_{disc}(\Delta)$ can simply be defined as all the isolated eigenvalues\footnote{i.e., $\Delta\phi = \lambda\phi$ for some $\phi\in\text{Dom}(\Delta)$ and $(\lambda -\epsilon, \lambda +\epsilon)\cap\sigma(\Delta) = \lambda$.} in $\sigma(\Delta)$ with finite multiplicity and the essential spectrum $\sigma_{ess}(\Delta)$ is simply the complement.  The essential spectrum is stable under compact perturbations in the sense that 
$\sigma_{ess}(\Delta + T) = \sigma_{ess}(\Delta)$ for any compact operator $T$.  One would like to think that compact perturbations of the metric on a manifold amounts to perturbing $\Delta$ by a compact operator.  This is a useful idealization.  However, it is heuristic.     

\quad In light of the heuristic notion that the essential spectrum is invariant under compact perturbation of the metric, the following result (for which its less general version first appeared in \cite{dek}) is not surprising.

\begin{theorem}[see \cite{dek} or ~\cite{LL-1}]\label{lowerbd}
Let $\Omega$ be a quantum layer over an isometrically immersed surface in $\mathbb{R}^3$, if we further assume that $\|\vec{A}\|\to 0$ at infinity, then the bottom of the essential spectrum of the Dirichlet Laplaician is bounded below as 
\[
 \inf\sigma_{ess}(\Delta) \geq \Big(\frac{\pi}{2a}\Big)^2.
\] 
\end{theorem}

The variation principle says that 
\[
\inf\sigma(\Delta) = \inf_{\phi\in C_0^{\infty}(\Omega)}
\frac{\int_{\Omega}|\nabla\phi|^2}{\int_{\Omega}\phi^2}.
\]
Therefore assuming asymptotic flatness of the surface $\varSigma$, to prove that the quantum layer $\Omega$ has ground state it suffices to find a test function so that
\begin{equation}\label{above} 
\inf\sigma(\Delta) < (\pi/2a)^2.
\end{equation}  
This was achieved first in \cite{dek}, and more generally in \cite{LL-1}(Theorem 1.1) by assuming $L^1$ Gauss curvature and non-positivity of the total Gauss curvature.  Thus there is the following conjecture: 
\newtheorem{conjecture}{Conjecture}
\begin{conjecture}\label{conj1}
Let $\varSigma$ be a complete, non-compact surface isometrically immersed in $\mathbb{R}^3$ and asymptotically flat.  Consider the quantum layer built over it.  If $\varSigma$ is not the plane and the Gaussian curvature $K$ is integrable on $\varSigma$, does the Dirichlet Laplacian on the quantum layer have non-empty discrete spectrum?  
\end{conjecture}
The remaining case of the conjecture above, which is the $\int_{\varSigma}K > 0$ case, was answered partially through special examples of surfaces in both \cite{dek} and \cite{LL-1}.  In particular, in \cite{LL-1}(Theorem 1.3) the example of quantum layers over a convex surface (graph of a convex function $f:\mathbb{R}^2 \longrightarrow \mathbb{R}$ in $\mathbb{R}^3$) is demonstrated to have ground state (assuming asymptotic flatness).  
A convex surface has positive Gauss curvature everywhere.  

\quad We must emphasize that in the proof of (\ref{above}) for the $\int_{\varSigma}K \leq 0$ case, there was no reference to the immersion of $\varSigma$ in the ambient Euclidean space.  Thus the existence of ground state there is of an intrinsic nature.  However, for the convex surface example, a careful analyis of the (integral of) mean curvature was used in an essential way.  Therefore, we believe that in order to answer the conjecture for the $\int_{\varSigma}K > 0$ case, the mean curvature - or more generally - the second fundamental form of $\varSigma$, must always play a central role.    

\quad In this paper we will discuss another rather general example of surfaces whose layers possess bound states, and the mean curvature on such surfaces will be heavily involved in the analysis.  We consider embedded surfaces in $\mathbb{R}^3$ that is ruled outside a compact subset.  A ruled surface is such that for each point on the surface there passes an Euclidean straight line, called a ruling, lying also on the surface, and such that a local collection of such lines and the orthogonal flow lines through them constitute a local coordinate system on the surface.  Since this description is local in nature (at least we can allow the situation where the rulings end somewhere on the surface), we can define surfaces which are potentially only ruled outside a compact subset of the surface.  More details about ruled surfaces will be given in the next section.      
 
\quad Our main reslut is as follows:

\begin{theorem}[Main Theorem]\label{Main}
Let $\varSigma$ be a surface in $\mathbb{R}^3$ that is ruled outside a compact subset.  Then the bottom of the spectrum of the Dirichlet Laplacian on a quantum layer $\Omega$ of thickness $a$ over $\varSigma$ has the upperbound
\[
 \inf\sigma(\Delta) < \big(\frac{\pi}{2a}\big)^2.
\]
\end{theorem}

\newtheorem{remark}{Remark}
\begin{remark}
Note that for a surface which is ruled outside a compact set, the Gauss curvature is automatically integrable. The surface also can't be flat outside a compact set, because otherwise the total Gauss curvature would be zero and in that case, the theorem follows from the main theorem in~\cite{CEK-1}, or ~\cite{LL-1}.
\end{remark}

Using Theorem \ref{lowerbd}, we also obtain the following:

\newtheorem{corollary}{Corollary}\label{cor1}
\begin{corollary}
The layer above, along with the assumption that the second fundamental form goes to zero at infinity, has Dirichlet ground state.  
\end{corollary}

\quad We end this section with an overview of the rest of the paper.  In section \ref{grs} we discuss the (local) geometry of ruled surfacs.  In section \ref{topology} we give relevant information about the topology of noncompact surfaces with integrable Gauss curvature.  The dicussion there will center around the generalized Gauss-Bonnet Theorem of Hartman \cite{hartman}.  In particular, we will make essential use of the theorem of B. White \cite{bwhite} in the proof of the main theorem.  Section \ref{proof} contains the proof of the 
main theorem.  

{\bf Acknowledgement.}  The first author would like to thank his program colleagues Spiros Karigiannis and William Wylie for enlightening discussions on surface topology during his stay at MSRI for the 2006-2007 year.  

\section{The Geometry of Ruled Surfaces}\label{grs}
\quad Here we discuss some basic geometry about ruled surfaces.

\begin{definition}\label{def1}
A non-intersecting surface $\varSigma$ in $\mathbb{R}^3$ is called a ruled surface if every point lies in a coordinate chart of the form 
\begin{equation}\label{ruled}
 x(s,v) = \beta(s) + v\delta(s),  
\end{equation}
where $\beta$ and $\delta$ are vectors in $\mathbb{R}^3$.
\end{definition}

It is important to note that the coordinate charts described above are in general only local.   We can always choose $\beta$ to be unit speed and $\delta$ to be a unit vector field.  Furthermore, we may reparameterize the coordinate chart above so that
\begin{equation}\label{special}
\begin{cases}
 |\beta'| = 1;\\
 |\delta| = 1;\\
 \langle\beta', \delta\rangle = 0.
\end{cases}
\end{equation}

By product rule we also have $\langle\delta^{'}, \delta\rangle = 0$.  From Definition \ref{def1} we see that at each point of a ruled surface $\varSigma$, there exists a straight line segment in $\mathbb{R}^3$ that also lies in $\varSigma$.  Such a line segment, called a generator, is simply the $v$-parameter curve of the coordinate system provided by (\ref{ruled}).  Now, we see that 
$\langle x_s,x_v\rangle = \langle\beta^{'},\delta\rangle + v\langle\delta^{'},\delta\rangle = 0$, hence we have local orthogonal coordinate systems on a ruled surface $\varSigma$.  Then letting $X = x_v$ and $Y = x_{s}/|x_{s}|$, we get local orthonormal frames on $\varSigma$. Denote by $f(v)$ and $h(t)$ the integral curves of $X$ and $Y$, respectively.  Let $N$ denote the oriented unit normal on $\varSigma$, then the Gauss curvature on $\varSigma$ is 
\begin{equation}\label{gauss}
K = \langle\frac{d}{dv}N\circ f,\, X\rangle\langle\frac{d}{dt}N\circ h,\, Y\rangle \, - \, 
\langle\frac{d}{dv}N\circ f,\, Y\rangle^2.
\end{equation}  
Since the ruling $f(v)$ is really a line segment, $f^{''} = 0$, so by product rule we see that
\[
 \langle\frac{d}{dv}N\circ f,\, X\rangle = -\langle N, f^{''}\rangle = 0. 
\]
Thus $K\leq 0$ on a ruled surface by (\ref{gauss}).  A surface $\varSigma$ is called a \text{\it{developable surface}} if it is a ruled surface such that its normal is parallel in $\mathbb{R}^3$ along any of its generators, i.e. $\frac{d}{dv}N\circ f = 0$.  Then again by (\ref{gauss}) we see that $K\equiv 0$ on a developable surface.  In fact the concept of a  developable surface is inseparable from the zero Gauss curvature condition:

\begin{theorem}[Massey's Theorem]\label{Massey}\, (Corollary to Theorem 5, \cite{hicks}) \quad A complete, connected surfcace in $\mathbb{R}^3$ is a developable surface if and only if $K\equiv 0$.\footnote{It is an elementary exercise to show that a ruled surface is developable if and only if it has zero Gauss curvature.  The more general theorem of Massey, on the other hand, is quite non-trivial.}
\end{theorem}

Next we compute the mean curvature on a ruled surface $\varSigma$.  Let $\times$ denote the cross pruduct of vectors in $\mathbb{R}^{3}$.  Moreover, since $|x_v| = |\delta| = 1$ we have $|x_s \times x_v| = |x_s|$.  The mean curvature on $\varSigma$ is computed as 
\begin{align}\label{mean}
H &= \langle-\frac{d}{dv}N\circ f,\, X\rangle + \langle-\frac{d}{dt}N\circ h,\, Y\rangle\notag\\
    &= \frac{\langle x_s \times x_v , \, x_{ss}\rangle}{|x_s\times x_v|^3}.
\end{align}

\section{The Topology of Noncompact Surfaces}\label{topology}
\quad The famous Gauss-Bonnet Theorem asserts that the total Gauss curvature of a compact $2$-dimensional manifold without boundary is a constant multiple of its Euler Characteristic.  If a surface is complete, non-compact, but has integrable Gauss curvature, then the total Gauss curvature is no longer completely topological.  In 1957, A. Huber showed that a complete, noncompact surface $\varSigma$ of integrable curvature is conformally equivalent to a compact Riemann surface with finitely many punctures.  In particular, Huber also showed that in this case
\begin{equation}\label{Huber}
 \int_{\varSigma}K \leq 2\pi\chi(\varSigma),
\end{equation}
where $\chi(\varSigma)$ is the Euler characteristic of $\varSigma$.  The finitely many punctures correspond to the ends of $\varSigma$.  Let us denote the ends by $\{E_1,...,E_k\}$, and define the corresponding isoperimetric constants 
\[
 \lambda_i = \lim_{r\to\infty}\frac{\text{Area}(B(r)\cap E_i)}{\pi r^2},
\]
relative to any fixed point $p\in\varSigma$ with respect to which the geodesic distance $r$ 
is measured.  The ends contribute to the deficit in (\ref{Huber}) via the following wonderful result.

\begin{theorem}[Hartman \cite{hartman} '64]
Let $\varSigma$ be a complete, noncompact surface with integrable Gauss curvature.  Then 
\begin{equation}\label{Hartman}
\frac{1}{2\pi}\int_{\varSigma}K = \chi(\varSigma) - \sum_{i=1}^{k}\lambda_i, 
\end{equation}
where $\chi(\varSigma)$ is the Euler characteristic of the surface.
\end{theorem}

The Euler Characteristic is defined by $\chi(\varSigma) = \sum (-1)^{i}b_i$, where $b_i$ is the $i$-th Betti number.  Now, we always assume implicitly that manifolds are connected, hence path-connected.  Thus for a noncompact (path-connected) surface, $\chi(\varSigma) = 1 - b_1$, and in the case that $\varSigma$ also has integrable Gauss curvature we see that in fact 
\begin{equation}
 \int_{\varSigma}K \leq 2\pi.
\end{equation}
Moreover, if we assume that $\int_{\varSigma}K >0$, then by (\ref{Hartman}) we must have 
$b_1 =0$ as well.  This means $\chi(\varSigma) = 1$, which via the uniformization theorem for surfaces implies that $\varSigma$ is conformally equivalent to $\mathbb{R}^2$.  Therefore we see that positve total Gauss curvature surfaces are topologically very simple (it is just $\mathbb{R}^2$).  However, ironically the analysis which are hopeful towards deducing the existence of ground state is significantly less straight-forward.  Contrasting this with the non-positive total Gaussian curvature case, where via Hartman's formula, if we start with any surface and add sufficient many handles we could obtain a ground state eventually.\footnote{This was first discussed in a subsequent paper  \cite{dek}  to \cite{CEK-1}.}

\quad For later use, we shall need the following result of B. White, which depends on the way in which a surface sits in $\mathbb{R}^3$.
\begin{theorem}[B. White \cite{bwhite} '87]\label{bbwhite}
Let $\varSigma$ be a surface immersed in $\mathbb{R}^3$.  If $\int_{\varSigma}|\vec{A}|^2 < \infty$, then $K$ is integrable and $\int_{\varSigma}K = 4\pi n$ for some $n\in\mathbb{Z}$.
\end{theorem}

\section{Proof of Main Result}\label{proof}
\quad   
\quad We assume that the surface $\varSigma$ is ruled outside $\overline{B(R_0)}$, where we have suppressed the reference point $x_0$ with respect to which the geodesic distance $R_0$ is measured.  From now on if the center of a geodesic ball of certain radius is suppressed, it is implicit that the center is the point $x_0$.  Thus each point in $\varSigma \setminus \overline{B(R_0)}$ is contained in a local coordinate chart given by $x(s,v) = \beta(s) + v\delta(s)$ for a curve $\beta$ in $\mathbb{R}^3$ and a nonzero outward-pointing vector field $\delta$ along $\beta$.   

\quad We would like to have a finite cover of $\varSigma \setminus \overline{B(R_0)}$ with ruled coordinate charts satisfying property (\ref{special}).  Let $R_1 > R_0$.  For any point $p\in\partial B(R_1)$, there exists a local ruled coordinate chart with property (\ref{special}).  However, we must pay attention to the possibility that in the reparameterization described above, 
$t(s)$ may be so negative for some $s$ that $\gamma(s)$ is no longer 
on the surface $\varSigma$.  Such a possibility requires consideration since our surface is not entirely ruled.  Note that however, 
\begin{equation}\label{control}
 |t(s)| \leq \int_0^s \big|\langle\beta',\delta\rangle\big|\,du \leq \int_{0}^s\,du = l(\beta). 
\end{equation}  
Therefore for each point $p\in\partial B(R_1)$, without loss of generality let 
$x(s,v): (-\epsilon, \epsilon)\times(b,\infty) \longrightarrow \varSigma$ be a ruled coordinate chart satisfying property (\ref{special}) such that $p=x(0,0)$.  If we further require that $\epsilon < l(\beta)$, then the image of $x$ will definitely lie in $\varSigma$.\footnote{The number $b < 0$ can always be chosen accordingly so that this is true.}  Then since 
$\partial B(R_1)$ is compact, there is a finite collection of such coordinate charts covering $\varSigma\setminus\overline{B(R_1)}$, assuming that we also take $R_1$ large enough.    

\quad By (\ref{mean}) and using property (\ref{special}), the mean curvature has the expression
\begin{align}\label{meancurv}   
 H &= \frac{\langle x_s \times x_v , \, x_{ss}\rangle}{|x_s\times x_v|^3}\notag\\
   &=\frac{\langle\beta^{'}\times\delta,\beta^{''}\rangle + v\Big(\langle\delta^{'}\times\delta,\beta^{''}\rangle + \langle\beta^{'}\times\delta,\delta^{''}\rangle\Big) + v^2\langle\delta^{'}\times\delta,\delta^{''}\rangle}{\Big(1+2v\langle\beta^{'},\delta^{'}\rangle + v^2|\delta^{'}|^2\Big)^{3/2}}.
\end{align}
Let us denote the numerator of $H$ above by $P$, and the term inside the $3/2$ power in the denominator by $Q$.  We will simplify matters by using the notations:
\begin{align*}
& A=2\langle\beta',\delta'\rangle\,\,;\,\,B=|\delta'|^2\,\,;\,\,
 D=\langle\beta'\times\delta,\beta^{''}\rangle\,\,;\,\,\\&
 E=\langle\delta'\times\delta,\beta^{''}\rangle+\langle\beta'\times\delta,\delta^{''}\rangle
 \,\,;\,\,J=\langle\delta'\times\delta,\delta^{''}\rangle.
\end{align*}

 Note that since we assume $\varSigma$ is non-flat on $\varSigma\setminus \overline{B(R_0)}$, there must be a point $p\in\varSigma\setminus \overline{B(R_0)}$ such that $H(p)\ne0$.  Moreover, by (\ref{meancurv}) if we assume $v$ is large enough, $H \ne 0$ along the ruling line for $s=0$.  
Switching the orientation of $\varSigma$ if necessarily, we can then take $\epsilon$ small enough and $v_0$ large enough for all $s\in(-\epsilon,\epsilon)$ such that $H>0$ for all 
$v\geq v_0$ on 
$x\big((-\epsilon,\epsilon)\times(v_0,\infty)\big)$.  This assumption on the positive sign of $H$ will be made implicitly throughout the remaining arguments.

\newtheorem{lemma}{Lemma}
\begin{lemma}\label{lem1}
 Let $t > 0$ and $\alpha, \beta\in(-\epsilon,\epsilon)$ such that $\alpha < \beta$.  Suppose that deg($P$) and deg($Q$) remain unchanged on $(\alpha,\beta)$.  Then the integral
 \[
 \int_{\alpha}^{\beta}\int_{t}^{t+t_0}H |x_s\times x_v|\,dv\,ds = o(t_0) 
 \]
 for large enough $t$ if and only if deg($P$) $<$ deg($Q$) for all $s\in(\alpha,\beta)$. 
\end{lemma}

\text{\bf{Proof.}}
\quad First we suppose deg($P$) $<$ deg($Q$), which consists of only two cases: 
$\big($deg($P$),\,deg($Q$)$\big)$ = $(1,2)$, and $(0,2)$.  For the first case, since 
$[\alpha, \beta]$ is compact we can choose constants $C_1, C_2>0$ such that for large enough $v$, 
\[
 P < C_1 v \hskip 1cm \text{and} \hskip 1cm  Q > C_2 v^2 
\]
for all $s\in(-\epsilon,\epsilon)$.  Then by (\ref{meancurv}), 
since $|x_s\times x_v| = \sqrt{1+Av+Bv^2}$, we see that 
\begin{equation}
\int_{\alpha}^{\beta}\int_{t}^{t+t_0}H |x_s\times x_v|\,dv\,ds \,= \, \int_{\alpha}^{\beta}\int_{t}^{t+t_0}\frac{P}{Q}\,dv\,ds 
\,< \, C\log\Big(\frac{t+t_0}{t}\Big)
\end{equation} 
for some constant $C>0$, hence $o(t_0)$.  For the second case, similarly we can find constants 
$C_1, C_2 >0$ such that for large enough $v$,
\[
 P < C_1 \hskip 1cm \text{and} \hskip 1cm Q> C_2 v^2
\]
for all $s\in(\alpha,\beta)$.  Then we see that for large enough $t$,
\begin{equation}
\int_{\alpha}^{\beta}\int_{t}^{t+t_0}H |x_s\times x_v|\,dv\,ds 
\,< \, C\Big(\frac{1}{t} - \frac{1}{t+t_0}\Big),
\end{equation} 
which is clearly $o(t_0)$. 

\quad Next suppose deg($P$) $\geq$ deg($Q$), which consists of three cases:
 $\big($deg($P$),\,deg($Q$)$\big)$ = $(2,2)$, $(0,0)$, and $(1,0)$.  In any of these cases, we can long divide (for each $s$) and get 
\begin{equation}\label{longdivide}
 \frac{P}{Q} = \tilde{P} + \frac{R}{Q},
\end{equation}
where deg($\tilde{P}$) = deg($P$) $-$ deg($Q$), and deg($R$) $<$ deg($Q$).  For large enough $t$, the integral over $(\alpha, \beta)\times (t, t+t_0)$ of the second term in (\ref{longdivide}) is $o(t_0)$ by applying the previous argument.  However, integrating the first term in (\ref{longdivide}) we see that it becomes a polynomial in $t_0$ of degree at least $1$.  To be more precise, for large eough $t$ we have
\[
 \int_{\alpha}^{\beta}\int_{t}^{t+t_0}H |x_s\times x_v|\,dv\,ds = f(t_0) + o(t_0),
\]   
where $f(t_0)$ is a polynomial of degree either $1$ (corresponding to the (2,2) and (0,0) cases listed above) or $2$ (corresponding to the (1,0) case) with coefficients in $t$.  Hence it cannot be $o(t_0)$.   
\qed  

\quad   Lemma \ref{lem1} is an important lemma that will allow us to deduce the existence of ground states.  In particular, we will need the integral of $H$ over a sector-like region in $\varSigma\backslash\overline{B(R_0)}$ to grow at least linearly in $t_0$ (i.e., in the ruling direction).  In view of this, we shall need to prove the following result. 
\newtheorem{proposition}{Proposition}
\begin{proposition}\label{daprop}
Let $\varSigma$ be a surface embedded in $\mathbb{R}^3$ that is ruled, but non-flat outside $\overline{B(R_0)}$ for some $R_0 > 0$.  Furthermore, assume the total Gauss curvature $\int_{\varSigma}K \,>\, 0$.  Then there exists $p\in\varSigma\setminus \overline{B(R_0)}$, $H(p)\ne0$, and a ruled coordinate chart 
$x: (-\epsilon, \epsilon)\times(b,\infty)\longrightarrow \varSigma$ satifying (\ref{special}) such that $p = x(0,0)$ and deg($P$) $\geq$ deg($Q$) at $s =0\in(-\epsilon,\epsilon)$.  
\end{proposition}
\quad Before we proceed with the proof of Proposition \ref{daprop}, let us remark on its nature.  Any point 
$p\in\varSigma\setminus \overline{B(R_0)}$ such that $H(p)\ne0$ lies in a ruled coordinate chart satifying (\ref{special}), and we can always reparameterize the chart so that $x(0,0) = p$ by translation in the $v$-direction.  Moreover, the degrees of $P$ and $Q$ certainly do not change under translation in $v$.  Therefore, the proposition is really a statement about the ruled coordinate charts covering $\varSigma\setminus \overline{B(R_0)}$.   

\text{\bf{Proof of Proposition \ref{daprop}.}} 
\quad Let $p\in\varSigma\setminus \overline{B(R_0)}$ with $H(p)\ne0$.  Consider a ruled coordinate chart $x: (-\epsilon, \epsilon)\times(b,\infty)\longrightarrow \varSigma$ containing $p$ such that $x(0,0) = p$.  Suppose deg($P$) $<$ deg($Q$) for all $s\in(-\epsilon,\epsilon)$ on such a chart, we will prove that 
\begin{equation}\label{dakey}
\int_{-\epsilon}^{\epsilon}\int_{0}^{\infty} H^2 |x_s\times x_v|\,dv\,ds \, < \, \infty.
\end{equation}

\quad Now suppose deg($P$) $<$ deg($Q$) at $s=0$ for all points $p=x(0,0)$ with $H(p)\ne0$.  By continuity, we may assume $\epsilon$ is small enough so that $H(s,0)\ne 0$ for all $s\in(-\epsilon,\epsilon)$.  For such a chart $x$, deg($P$) $<$ deg($Q$) for all $s\in(-\epsilon,\epsilon)$.  To see this, suppose 
deg($P$) $\geq$ deg($Q$) at some $s_0\in(-\epsilon,\epsilon)$.  Then we can make the reparametrization 
\[
 \beta(s) \longrightarrow \beta(s-s_0) \hskip 1cm  \text{and} \hskip 1cm 
 \delta(s) \longrightarrow \delta(s-s_0)
\]
to obtain a new ruled coordinate chart $\tilde{x}$ with $\tilde{x}(0,0) = x(s_0,0)$.  Note 
that property (\ref{special}) is preserved for $\tilde{x}$, and moreover 
the coefficients $A, B, D, E, J$ are also invariant under the reparametrization above.   
Thus deg($\tilde{P}$) = deg($P$) and deg($\tilde{Q}$)=deg($Q$), and we would get a contradiction.  Then via (\ref{dakey}) and by our earlier remark on the existence of a finite cover of $\varSigma\setminus \overline{B(R_0)}$ by ruled coordinate charts, it would imply that $H\in L^2(\varSigma)$.  By the elementary formula $H^2 = |\vec{A}|^2 + 2K$ and the fact that $K\in L^1(\varSigma)$, we must then have $\int_{\varSigma}|\vec{A}|^2 \,<\, \infty$, which by B. White's theorem (Theorem~\ref{bbwhite}) means that $\int_{\varSigma}K = 4\pi n$ for some $n\in\mathbb{Z}$.  However, the assumption of connectedness and noncompactness of $\varSigma$ implies we must have Euler characteristic $\chi(\varSigma) \leq 1$.  Then by Hartman's formula and our hypothesis we must have  
\begin{equation}\label{realisation}
 0\,<\,\int_{\varSigma}K\,\leq\,2\pi, 
\end{equation}     
which would give a contradiction.    

\quad Now we prove (\ref{dakey}) assuming deg($P$) $<$ deg($Q$) for all 
$s\in(-\epsilon,\epsilon)$ on a chart.  First we see that
\begin{equation}
 \int_{0}^{t} H^2|x_s\times x_v|\,dv = \int_{0}^{t} \frac{P^2}{Q^{5/2}}\,dv.
\end{equation}
First consider all $s\in(-\epsilon,\epsilon)$ such that 
$\big($deg($P$),\,deg($Q$)$\big)$ = $(1,2)$.  Since the coefficients of $P$ and $Q$ are functions of $s$ bounded on $(-\epsilon,\epsilon)$, we can choose $c>0$ such that $Q>1$, $Q^2 > C v^4$ for some 
constant $C > 0$, and 
$P^2 < C v^2$ for some constant $C > 0$, for all $v>c$ and for all such $s$.  Then for any such $s$, we have  
\begin{align}\label{1,2}
 \int_{c}^{t}\frac{P^2}{Q^{5/2}}\,dv \,&\leq\, \int_{c}^{t}\frac{P^2}{Q^{2}}\,dv\notag\\
 &< C\int_{c}^{t}\frac{1}{v^2}\,dv\notag\\
 &= C\Big(\frac{1}{c} - \frac{1}{t}\Big),
\end{align}  
which converges as $t\to\infty$.  Similarly, for any $s\in(-\epsilon,\epsilon)$ at which $\big($deg($P$),\,deg($Q$)$\big)$ = $(0,2)$ we have
\begin{align}\label{0,2}
\int_{c}^{t}\frac{P^2}{Q^{5/2}}\,dv \,&\leq\, P^2\int_{c}^{t}\frac{1}{Q^{2}}\,dv\notag\\
 &< C\,P^2\int_{c}^{t}\frac{1}{v^4}\,dv\notag\\
 &= C\,P^2\Big(\frac{1}{c^3} - \frac{1}{t^3}\Big),
\end{align}
which again converges as $t\to\infty$.  Combining (\ref{1,2}) and (\ref{0,2}) proves (\ref{dakey}) and hence the proposition.
\qed

\begin{remark}
\quad Although this is not necessary for the proof of Proposition \ref{daprop}, it is in good spirit to check that deg($P$) $\geq$ deg($Q$) at $p=x(0,0)$ does not render $H^2$ unintegrable over $(-\epsilon,\epsilon)\times(0,\infty)$ for a small enough $\epsilon >0$ where 
deg($P$) and deg($Q$) are unchanged.  As in the previous case of deg($P$) $<$ deg($Q$), we can narrow down to three cases.  The first two consist of deg($P$) $= 1,0$ coupled with deg($Q$) $= 0$.  It is clear in these two cases that
\[ 
\int_{-\epsilon}^{\epsilon}\int_{0}^{t}H^2|x_s\times x_v|\,dv\,ds
\]
does not converge as $t\to\infty$.  For the last case of deg($P$) $=$ deg($Q$) $= 2$, one can verify that there exists $c,k > 0$ such that 
\begin{align}
 \int_{c}^{t}\frac{P^2}{Q^{5/2}}\,dv &\geq k\int_{c}^{t}\frac{1}{v}\,dv \notag\\ 
 &= k\big(\log t - \log c\big) \longrightarrow \infty
\end{align} 
as $t\to\infty$.
\end{remark}

\begin{corollary}\label{dacor}
There exists a point $p\in\varSigma\setminus\overline{B(R_0)}$, $H(p)\ne 0$, and a ruled coordinate chart 
$x: (-\epsilon,\epsilon)\times(b,\infty) \longrightarrow \varSigma$
satisfying (\ref{special}) such that deg($P$) $\geq$ deg($Q$) for all $s\in(-\epsilon,\epsilon)$.  Moreover, we can choose the chart so that deg($P$) and deg($Q$) are fixed for all $s\in(-\epsilon,\epsilon)$.
\end{corollary}

\text{\bf{Proof.}}\quad 
By Proposition \ref{daprop}, there exists at least one point $p\in\varSigma\setminus\overline{B(R_0)}$, $H(p)\ne 0$, and a ruled coordinate chart $x:(-\epsilon,\epsilon)\longrightarrow\varSigma$ such that $x(0,0)=p$ and 
deg($P$) $\geq$ deg($Q$) at $s=0$.  Suppose at such a point $p = x(0,0)$, the first assertion of the corollary is false.  Then by the smoothness of the coefficients of $P$ and $Q$ in $s$, there must exist an $\overline{\epsilon}>0$, $\overline{\epsilon}\leq\epsilon$, such that   
deg($P$) $\geq$ deg($Q$) at $s=0$ and deg($P$) $<$ deg($Q$) for all 
$s\in(-\overline{\epsilon},\overline{\epsilon})\setminus\{0\}$.  Thus the ruling lines that pass through these points must comprise a set of discrete lines, and hence of measure zero.  Then applying the integration of $H^2$ argument in Proposition \ref{daprop}, we would get the same contradiction as we did there.

\quad Next we argue that we can fix deg($P$) and deg($Q$) in a small enough interval of $s$.  Observe that due to the smoothness of the coefficient functions in $s$, the degrees of $P$ and $Q$ cannot decrease in an arbitrarily small neighborhood of $s$, but it can certainly increase.  Now, if $\big($deg($P$),\,deg($Q$)$\big)$ $= (2,2)$ at $s=0$, then since this is the case of the largest possible degrees, we can certainly find an $\epsilon > 0$ small enough so that the degrees remain $2$ for all $s\in(-\epsilon,\epsilon)$.  If we are in the $(1,0)$ case at $s=0$, then either the degrees remain as $(1,0)$ in a small interval about $s=0$, or there exists an $s_0$ near $s$ at which the degrees increase to $(2,2)$ and to which we can apply the preceding argument upon reparametrization of the chart.  For the last case of $(0,0)$, if the degrees increase near $s=0$ then we simply apply the arguments in the previous two cases.  
\qed

\text{\bf{Proof of Theorem \ref{Main}.}}
\quad Now, the surface is not totally-geodesic (non-planar) by our hypothesis.  Therefore if $\int_{\varSigma}K \leq 0$, the conclusion of the theorem follows as a special case of the result in \cite{LL-1}.  For the remaining of the proof we will assume that $\int_{\varSigma}K > 0$.

\quad By the variation principle,it suffices to find a test function 
$\phi\in W_{0}^{1,2}(\Omega)$ such that 
\begin{equation}\label{var}
Q(\phi,\phi) = \int_{\Omega}|\nabla\phi|^2 -\Big(\frac{\pi}{2a}\Big)^2\int_{\Omega}\phi^2 
< 0.
\end{equation}
We define a test function (family of test functions) of the form\footnote{This form of test function first appeared in \cite{dek}, and was also used extensively in \cite{LL-1}.  However, the argument we will give here for (\ref{var}) is different in an essential way.}  
\begin{equation}\label{test}
 \phi_{\epsilon} = \chi\psi + \epsilon\chi_1 j,
\end{equation}
where $\epsilon$ is some non-zero number to be determined, $\chi = \cos\frac{\pi}{2a}u$, $\chi_1 = u\cos\frac{\pi}{2a}u$, and $\psi$ and $j$ are defined below.

\quad The assumption of integrable Gauss curvature implies that $\varSigma$ is parabolic\footnote{For a proof of this, see \cite{LL-1}.}, and hence for any $R_1 > 0$ there exists an $R_2 > R_1$ for which 
\[
 \int_{\varSigma}|\nabla\psi_{R_1, R_2}|^2 \, < \, \frac{\epsilon_0}{2},
\] 
where $\psi_{R_1, R_2}$ is the unique solution to the boundary value problem
 \begin{equation}\label{psi1}
  \begin{cases}
   \Delta \psi = 0 \quad & \text{on $B(R_2)\setminus B(R_1)$};\\
   \psi|_{B_p(R_1)} \equiv 1; \\
   \psi|_{\varSigma\setminus B_p(R_2)} \equiv 0.
  \end{cases}
 \end{equation}
We will let $R_1 > R_0$.  Then for $R_1 < R_2 < R_3 < R_4$, we let 
\[
 \psi = \psi_{R_3,R_4} - \psi_{R_1,R_2}.
\]

\quad We want $j$ to be a $W^{1,2}$ function on $\varSigma$ with support in $\{\psi\equiv 1\}$ and $j\leq 1$.  Before defining $j$ precisely, we proceed with some preliminary estimates.  By our choices of $\chi$ and $\chi_1$, the fact that $K\leq 0$ on $\varSigma\setminus B(R_0)$, and the requirement that $j\leq 1$ and $\text{supp}j \subset \{\psi\equiv 1\}$, we get 
\begin{align}\label{1stineq}
Q(\phi_\epsilon,\phi_\epsilon) &= Q(\chi\psi,\chi\psi) + 2\epsilon Q(\chi\psi,\chi_1 j) + 
\epsilon^2 Q(\chi_1 j, \chi_1 j)\notag\\
&\leq C_1 \int_{\varSigma}|\nabla\psi|^2 + \epsilon a\int_{\varSigma}jH + \epsilon^2 C_2\|j\|_{W^{1,2}}^2
\end{align}
for $C_1$, $C_2 \,>\,0$ depending only on the geometry of $\Omega$.  We will choose $j$ so that $\|j\|_{W^{1,2}} \ne 0$.  Then viewing the right-hand side of the inequality in (\ref{1stineq}) as a quadratic polynomial in the variable $\epsilon$, it will be negative for some $\epsilon$ if and only if its discriminant is positive, which is equivalent to the condition
\begin{equation}\label{cond}
\frac{\big(\int_{\varSigma}jH\big)^2}{\|j\|_{W^{1,2}}^2} \, > \, 
C_{1}\int_{\varSigma}|\nabla\psi|^2,
\end{equation}     
where we absorbed all the geometric constants into a single constant $C_1$.  Now, by our choice of $\psi$ along with the parabolicity of $\varSigma$, for $R_1 > R_0$ fixed we can choose $R_2$ and then $R_3 < R_4$ big enough so that 
\begin{align}\label{Q1}
 \int_{\varSigma}|\nabla\psi|^2 &= \int_{\varSigma}|\nabla\psi_{R_1,R_2}|^2 + \int_{\varSigma}|\nabla\psi_{R_3,R_4}|^2 \notag\\
 &< \frac{\epsilon_0}{2} + \frac{\epsilon_0}{2} = \epsilon_0. 
\end{align}
Observe that what is essential is the choice of $R_2$, as $R_3 < R_4$ can always be chosen after $R_2$ so that (\ref{Q1}) holds.  By our requirement on $j$, the choice of $R_2$ may affect the ratio on the left-hand side of (\ref{cond}).  In view of this consideration, we seek a (family of) $j$ such that 
\begin{equation}\label{suff}
 \frac{\big(\int_{\varSigma}jH\big)^2}{\|j\|_{W^{1,2}}^2} \, \geq \, C
\end{equation}  
for a constant $C$ independent of $R_2$, as long as $R_2$ is large enough.  Inequality (\ref{suff}) is a sufficient condition for (\ref{cond}) since we can then choose $R_2$ large enough for a small enough $\epsilon_0$ satisfying (\ref{Q1}) and 
\[
 C_1 \epsilon_0 < C.
\]    
In a nutshell, the proof will be complete if we construct a (family of) $j$ so that (\ref{suff}) holds for some constant $C$ independent of $R_2$, for $R_2$ large enough.  

\quad By Corollary \ref{dacor}, let $p\in\varSigma\setminus\overline{B(R_0)}$ be a point such that $H(p)\ne0$ and consider a ruled coordinate chart $x: (-\epsilon, \epsilon)\times(b,\infty)\longrightarrow \varSigma$ satifying (\ref{special}) such that $p = x(0,0)$, deg($P$) $\geq$ deg($Q$), with deg($P$)and deg($Q$) fixed for all $s =0\in(-\epsilon,\epsilon)$.  Moreover, we let $R_1 = \text{dist}(x_0,p)$.  For any $R_2 > R_1$ and $t_0 >1$ let 
\[
 \Gamma = \{(s,v)\in\mathbb{R}^2 | -\epsilon \leq s \leq \epsilon \,,\, v_0 \leq v \leq v_0 + t_0\}
\]
such that $x\big(\Gamma\big) \subset B_p(R_3)\setminus B_p(R_2)$.    
We define $j$ by $j = j_1(s)j_2(v)$, with cut-off functions 
\begin{equation}\label{j1}
 j_{1}(s) = 
 \begin{cases}
  1 \quad \quad \alpha -\epsilon < s < \epsilon -\alpha;\\
  0 \quad \quad |s|\geq \epsilon
 \end{cases}
\end{equation} 
\quad and
\begin{equation}\label{j2}
 j_{2}(v) = 
 \begin{cases}
  1 \quad \quad v_0 + \alpha < v < v_0 + t_0 -\alpha;\\
  0 \quad \quad v\leq v_0,\, v\geq v_0 + t_0,
 \end{cases}
\end{equation}
where $\alpha >0$ is a fixed small number, $|j^{'}_1(s)|\leq 1/\alpha$, and 
$|j^{'}_2(v)|\leq 1/\alpha$. 

\quad Now, by the definition of $j$ above and (\ref{meancurv}), we have
\begin{equation}\label{intH}
\int_{\varSigma}jH > \int_{\{j\equiv 1\}}H \notag
= \int_{\alpha-\epsilon}^{\epsilon-\alpha}
\int_{v_0 + \alpha}^{v_0 + t_0 -\alpha}
\frac{P}{Q}\,dv\,ds.
\end{equation}
By our choice of the ruled coordinate chart $x$ above, Lemma \ref{lem1} implies 
\begin{equation}\label{limH}
\int_{\varSigma}jH > f(t_0) + o(t_0)
\end{equation}
for large enough $v_0$, where $f(t_0)$ is a poynomial in $t_0$ of degree $n = 1$ or $2$ and has coefficients in $v_0$.  

\quad Next we wish to give an upperbound estimate for $\|j\|_{W^{1,2}}$.  First we see that 
\begin{align}\label{jW12}
 \|j\|_{W^{1,2}} &= \int_{\varSigma}j^2 + \int_{\varSigma}|\nabla j|^2 \notag\\
 &\leq \big(1 + \|\nabla j\|_{\infty}^2\big)\text{vol}(x(\Gamma)).
\end{align} 
By the metric on $\varSigma$ given by the chart (\ref{ruled}), we see that 
\begin{align}\label{delj}
 |\nabla j|^2 \,\,&= \,\, 
 j_{2}\,^2G^{ss}\Big|\frac{\partial j_1}{\partial s}\Big|^2 + 
 2j_{1}j_2G^{sv}\frac{\partial j_1}{\partial s}\frac{\partial j_2}{\partial v} + 
 j_{1}\,^2G^{vv}\Big|\frac{\partial j_2}{\partial v}\Big|^2\notag\\
 &\leq \Big(\frac{1}{1+2v\langle\beta^{'},\delta^{'}\rangle+v^2|\delta^{'}|^2}\Big)\frac{1}{\alpha^2} + \frac{1}{\alpha^2}\notag\\
\end{align}
If the polynomial   
$1 + 2v\langle\beta^{'},\delta^{'}\rangle + v^2|\delta^{'}|^2$ has degree $2$ for all 
$s\in(-\epsilon,\epsilon)$, then since its discriminant
\[
 4\langle\beta^{'},\delta^{'}\rangle^2 - 4|\delta^{'}|^2 \,\,\leq\,\, 4|\beta^{'}|^2|\delta^{'}|^2 -4|\delta^{'}|^2 = 0,
\]
it is always positive for $v$ large enough.  The other possibility is that it is identically equal to $1$.  In any case, we can choose a $v_0$ large enough so that for all $v\geq v_0$ 
\begin{equation}\label{absj2}
 |\nabla j|^2 < C_3
\end{equation}
 for all $s\in(-\epsilon,\epsilon)$, for some constant $C_3 >0$. 

 \quad Next we will estimate the volume growth of $x(\Gamma)$.  The volume form in the ruled coordinate system is 
 \[
 d\varSigma = \sqrt{1+Av+Bv^2}\,ds\,dv.
 \]
 There are two possibilities at $s = 0$, either $B=0$ or $B\ne 0$.  In the latter case, we can certainly assume that $\epsilon$ is small enough so that $B\ne 0$ for all $s\in(-\epsilon,\epsilon)$.  If $B=0$ at $s=0$, then by a similar argument as in Corollary \ref{dacor}, we can assume that $B=0$ for all $s\in(-\epsilon,\epsilon)$.    
 
 \quad Suppose $B\ne 0$ for all $s\in(-\epsilon,\epsilon)$. If $A^2 -4B < 0$ at $s=0$ we can always take $\epsilon$ small enough so that $A^2 -4B < 0$ for all $s\in(-\epsilon,\epsilon)$.  Assuming so, we have
\begin{align}\label{volestm}
 \text{vol}(x(\Gamma))
  &= \int_{-\epsilon}^{\epsilon}\int_{v_0}^{v_0 +t_0}\sqrt{1+Av+Bv^2}\,\, dv \,ds\notag \\
 &= \int_{-\epsilon}^{\epsilon}\int_{v_0}^{v_0 +t_0}\sqrt{B(v+\frac{A}{2B})^2 + 1-\frac{A^2}{4B}}\,\, dv\, ds\notag\\
 &= \int_{-\epsilon}^{\epsilon}\int_{v_0 + \frac{A}{2B}}^{v_0 +t_0 + \frac{A}{2B}}
 \sqrt{Bx^2 + 1-\frac{A^2}{4B}}\,\, dx \,ds\notag \\
 &= \int_{-\epsilon}^{\epsilon}\sqrt{\frac{4B-A^2}{4B}}\int_{v_0 + \frac{A}{2B}}^{v_0 +t_0 + \frac{A}{2B}}
 \sqrt{\frac{4B^2}{4B-A^2}x^2 + 1}\,\,  dx\,ds\notag \\
 &= \int_{-\epsilon}^{\epsilon}\frac{4B-A^2}{4B^{3/2}}
 \int_{\sqrt{\frac{4B^2}{4B-A^2}}\big(v_0 +\frac{A}{2B}\big)}
 ^{\sqrt{\frac{4B^2}{4B-A^2}}\big(v_0 +t_0 +\frac{A}{2B}\big)} \sqrt{y^2 +1}\,\, dy\,ds\notag \\
 &\leq \int_{-\epsilon}^{\epsilon}\frac{4B-A^2}{4B^{3/2}}
 \int_{\sqrt{\frac{4B^2}{4B-A^2}}\big(v_0 +\frac{A}{2B}\big)}
 ^{\sqrt{\frac{4B^2}{4B-A^2}}\big(v_0 +t_0 +\frac{A}{2B}\big)} (y + 1)\,\, dy\,ds\notag \\
 & = t_0 v_0 \int_{-\epsilon}^{\epsilon}\sqrt{B}\, ds + 
 t_0^2\int_{-\epsilon}^{\epsilon}\frac{\sqrt{B}}{2}\, ds + 
 t_0\int_{-\epsilon}^{\epsilon}\Big(\frac{A}{2\sqrt{B}} + \sqrt{\frac{4B-A^2}{4B}} \Big)\, ds. 
\end{align}
On the other hand, if $B\ne 0$ for all $s\in(-\epsilon,\epsilon)$ and $A^2 -4B = 0$ at $s=0$, using the same argument in Corollary \ref{dacor} we can assume that $A^2 -4B =0$ for all $s\in(-\epsilon,\epsilon)$.  Assuming so, we have
\begin{align}\label{volestm2}
\text{vol}(x(\Gamma)) &= 
\int_{-\epsilon}^{\epsilon}\int_{v_0}^{v_0 +t_0} \sqrt{1+Av+Bv^2}\,dv\,ds \notag\\
&= \int_{-\epsilon}^{\epsilon}\int_{v_0}^{v_0 +t_0} 
\sqrt{B}\Big(v+\frac{A}{2B}\Big)\,dv\,ds \notag\\
&= t_0 v_0 \int_{-\epsilon}^{\epsilon}\sqrt{B}\, ds + 
 t_0^2\int_{-\epsilon}^{\epsilon}\frac{\sqrt{B}}{2}\, ds + 
 t_0\int_{-\epsilon}^{\epsilon}\frac{A}{2\sqrt{B}}\, ds. 
\end{align} 

For brevity, we will use the following notations:\\
\begin{align*}
&C_4 = \int_{-\epsilon}^{\epsilon}\sqrt{B}\, ds, \, 
C_5=\int_{-\epsilon}^{\epsilon}\frac{\sqrt{B}}{2}\, ds,\,\\
  &C_6=\int_{-\epsilon}^{\epsilon}\Big(\frac{A}{2\sqrt{B}} + \sqrt{\frac{4B-A^2}{4B}} \Big)\, ds\,\,\,or \,\,\int_{-\epsilon}^{\epsilon}\frac{A}{2\sqrt{B}}\, ds.
  \end{align*}

If $B=0$ for al $s\in(-\epsilon,\epsilon)$, then 
\begin{equation}\label{volestm3}
\text{vol}(x(\Gamma)) = 2\epsilon t_0.
\end{equation}

\quad Now, by (\ref{limH}), with any fixed $t_0 > 0$ there exists a $v_0$ large enough so that 
\begin{equation}\label{Hbelow}
 \int_{\varSigma}jH > C_7 \, t_0^n
\end{equation} 
for some constant $C_7 > 0$ (which depends on $v_0$), with $n = 1$ for cases (\ref{volestm}) or
(\ref{volestm2}) and $n=1$ or $2$ for case (\ref{volestm3}) corresponding to deg($P$) 
$=0$ or $1$.

Renaming the squre of $C_7$ as itself, by (\ref{jW12}), (\ref{absj2}), and (\ref{Hbelow}) with a choice of a large enough $v_0$, we see that in either case (\ref{volestm}) or case (\ref{volestm2}),
\begin{align}\label{fineq}
 \frac{\big(\int_{\varSigma}jH\big)^2}{\|j\|_{W^{1,2}}^2} > 
 \frac{C_7 t_0^{\,2}}{C_3\Big(t_0 v_0 C_4 + t_0^2 C_5 + t_0 C_6\Big)}.
\end{align}
Note that once $v_0$ is fixed, the constants $C_3, C_4, C_5, C_6$, and $C_7$ depend only on the metric along the $s$-parameter curve $\beta(s): (-\epsilon,\epsilon) \longrightarrow \varSigma$, which is fixed from the start.  Then the right-hand side of (\ref{fineq}) either converges to $C_7 / C_5$ (when $n=1$) or goes to infinity (when $n=2$), as we take $t_0\to\infty$ (by letting $R_3 \to \infty$).  Therefore for a large enough $t_0$, there must be a constant $C > 0$ such that  
\begin{equation}\label{vind}
 \frac{\big(\int_{\varSigma}jH\big)^2}{\|j\|_{W^{1,2}}^2} > C.
\end{equation}
For the case of (\ref{volestm3}), the denominator of the right-hand side of (\ref{fineq}) will always be linear in $t_0$, while the numerator is either qudratic or to the $4$th power in $t_0$, hence (\ref{vind}) is readily achieved.  The proof is now complete.
\qed

An immediate consequence of Theorem \ref{Main} is the following result, which follows via Massey's Theorem.

\begin{corollary}\label{K=0}
Let $\varSigma$ be an embedded surface in $\mathbb{R}^3$ with zero Gauss curvature but non-flat outside a compact subset.  Then for a layer $\Omega$ over $\varSigma$, we have
\[
 \inf\sigma(\Delta) < \Big(\frac{\pi}{2a}\Big)^2.
\]
Moreover, if the second fundamental form $\vec{A}\to 0$ at infinity on $\varSigma$ then ground states exist.  
\end{corollary}

\begin{remark}
Our proof depends heavily on the fact that the surface is {\it ruled} outside a compact set. If we assume that the surface is asymptotically flat, then it is {\it asymptotically ruled} on any compact set. However, with only the asymptotic flatness, we don't know the behavior of the surface at infinity. Thus the analysis in this paper is not enough to prove Conjecture~\ref{conj1}. More careful estimates are needed but the motivation of the development of the techniques  in this paper is clearly justified.
\end{remark}
\def\cprime{$'$} \def\cprime{$'$} \def\cprime{$'$}

\end{document}